\documentclass[onesidetide,10pt]{article}
\usepackage{amssymb}

\newtheorem{thm}{Theorem}[section]

\newtheorem{lem}[thm]{Lemma}

\newtheorem{deff}[thm]{Definition}

\input xy
\xyoption{all}

\def \k {\kern 5 pt}
\def \n {\noindent}

\def \e {\eta}

\def \N {\nabla}

\title{\bf Divergence theorems in path space III: hypoelliptic 
diffusions and beyond}
\begin{document}
\date{}
\maketitle

\begin{center} \author{ Denis Bell\footnote{Research partially 
supported by
NSF grant DMS-0451194}
\\Department of Mathematics, University of North Florida
\\4567 St. Johns Bluff Road South,Jacksonville, FL 32224,
U. S. A.
\\email: dbell@unf.edu}

\end{center}

\kern 50 pt

\n {\bf Abstract}.  Let $x$ denote a diffusion process defined on a
 closed compact manifold.
In an earlier article,
 the author introduced a new approach to constructing admissible
vector fields on the associated space of paths, under the assumption of 
ellipticity of $x$.
  In
this article, this method  is  extended to yield similar results
for degenerate diffusion processes.  In particular, these results apply 
to
non-elliptic diffusions satisfying H\"ormander's condition.

\vfil\break

\section{Introduction}





\n Let $X_1, \dots, X_n$ and $V$ denote  smooth
  vector fields on  a closed compact manifold $M$.
We fix a point $o\in M$ and a positive time $T$ and consider
the Stratonovich stochastic differential equation
 (SDE)
$$ dx_t = \sum_{i =1}^n
  X_i(x_t)\circ dw_i +V(x_t)dt, \kern 5 pt
 t \in [0, T] \eqno (1.1) $$

\kern - 10 pt

$$x_0=o. \kern 140 pt $$
where $w= (w_1,\dots,w_n)$ is a Wiener process in ${\bf R}^n$.
Assume  that the vector $V(x)$ lies within the span
 of 
$X_1(x),\dots, X_n(x) $, for all $x\in M$.
 The solution process
$x$                                     
is a random variable taking values in the space  of paths
$$
C_o(M) = \big \{\sigma: [0,T]\mapsto M/ \sigma (0)=o\big\},
$$
an
 infinite-dimensional manifold with tangent bundle consisting of fibers
$$
T_\sigma C_o(M) = \big\{r:[0,T]\mapsto TM/ \kern 2 pt
r_0=0, \kern 2 pt r_t\in T_{\sigma _t}M
\kern 2 pt  \forall
t \in [0,T]\big\}.
$$
The law $\gamma$ of $x$, as a measure on  $C_o(M)$,
can then  be considered as a generalized version of Wiener measure on 
$C_0({\bf R}^n)$.
A major goal in stochastic analysis is to extend the
 rich body of results that have been developed for the Wiener measure to this
non-linear setting.

The {\it Cameron-Martin space}, i.e. the space of paths
$\{\sigma: [0,T]\mapsto {\bf R}^n,\k \sigma_0 = 0\}$
with finite energy
$$
\int_0^T ||\dot \sigma_t||^2dt
$$
provides a geometrical framework for the Wiener measure and plays a central role
in its analysis.
Therefore, in addressing the problem raised above, it is natural to seek an analogue
 of the Cameron-Martin space
for the measure $\gamma$.
A reasonable  candidate for such an analogue is the set of vector fields on the
space $C_o(M)$
 that admit an  ``integration by parts" formula of the type
 described in the following
                           
\begin{deff}  A vector field $\eta$
on  $C_{o}(M)$ is admissible (with respect to  $\gamma$) if there   
  exists an $L^1$ function $Div(\e)$ such that the relation\footnote{The integrals
in (1.2) will usually
be written as expectations in the sequel.}
$$\int_{C_{o}(M)}\e(\Phi)d\gamma  = \int_{C_{o}(M)}\Phi Div(\e)d\gamma  \eqno (1.2)
$$  holds for a dense class of  smooth
  functions $\Phi$ on $C_{o}(M)$.
\end{deff}

\k

  The construction of admissible vector fields  is an important problem
that has been studied by many authors in the last three decades. A breakthrough
 in the problem
was achieved by  Driver [6] in 1992,  following important
partial results by Bismut [5]. Driver proved that  parallel translation along $x$
 of Cameron-Martin paths in
 $T_oM$ produces admissible vector fields on $C_o(M)$. A fundamental innovation in
[6] is
the use of the rotation-invariance of the Wiener process. This property also plays a crucial role in
the present work.

The work of Bismut and Driver stimulated a great deal of activity in this area
and the problem is still
being widely studied (cf., e.g. Driver [7], Hsu [9] and [10], Enchev
 \& Stroock [ 8], Elworthy, Le Jan \& Li [7]).
 Much of this work
has dealt with the
{\it elliptic} case, where the vector fields  $X_1,\dots, X_n$ in (1.1) are assumed to
span $TM$
 at all points of $M$. 
In  [1], the author introduced  a new approach to the problem
of constructing admissible vector fields on path space,
 again
in the elliptic setting.
The purpose of the present article,  the third in a series of papers on this theme (cf.
[1] and [2]),
 is to extend this approach to the case of {\it degenerate} diffusions.

The
 central object of study in the author's approach is the {\it It\^o map}
 $g: w\mapsto x$ defined by equation (1.1).  This  is used
to lift the problem from the manifold $M$
 to ${\bf R}^n$, where  classical integration by parts theorems can be
applied\footnote{This
method had previously been employed by Malliavin in his
probabilistic approach to the  hypoellipticity problem [12]. }.
``Lifting" is defined as follows.
\begin{deff}
 {\it A process $r$ taking values in ${\bf R}^n$   is said to be a
lift of $\eta$ to $C_0({\bf R}^n)$
 (via the
It\^o map)
if the following diagram commutes}
\footnote{Since  $g$ is non-differentiable
in the
classical sense the derivative $dg$  must be interpreted in
the
extended sense of  Malliavin. As this type of regularity is now generally 
well-understood by stochastic analysts,
this point will not be emphasized in the paper (cf. e.g.  the monographs [3], 
[13],  [14], [15] for an
introduction to the
Malliavin calculus).}
$$\xymatrix {&TC_0({\bf R}^n)  \ar[r]^{dg}  &TC_o(M)\\
             &C_0({\bf R}^n)\ar[u]^r \ar[r]_g   &C_o(M) \ar[u]_\e\\
}
$$
\end{deff}   
\noindent  The idea in [1]  is to simultaneously
construct a vector field
$\eta$ on $C_o(M)$ and an {\it admissible} lift $r$ of $\eta$ to $C_0({\bf R}^n)$.
In particular, this
requires that $r$ take the form
$$
r_t = \int_0^tA(s)dw_s + \int_0^tB(s)ds
$$
where $A$ and $B$ are continuous adapted processes taking values in 
$so(n)$ (the space of skew-symmetric $n\times n$ matrices)
and ${\bf
R}^n$  respectively. Processes 
of this form
thus comprise the tangent bundle $TC_0(R^n)$ in the above diagram.

 For  test functions\footnote{For test functions we use the set of smooth
cylindrical functions on $C_o(M)$.} $\Phi$  on
$C_o(M)$, one then has
$$
E\big [(\e\Phi)(x)\big ] = E\big [r(\Phi\circ g)(w)\big]$$
$$ = E\big [\Phi \circ
g(w)
 Div (r)\big ] 
$$
$$
= E\Big [ \Phi(x)E\big [Div(r)/x\big ]\Big ].
$$
where $Div$ denotes the divergence operator in the classical Wiener space.
Thus  $\e$ is admissible  with divergence
$$
Div(\e)(x) = E\big [Div(r)/x\big ].
$$

  An important consequence of the  ellipticity assumption is the fact 
that every non-anticipating
vector field
on $C_o(M)$ can be written in the form
$$
\eta_t = \sum_{i=1}^nh_i(t)X_i(x_t)
$$
where $h_i, i=1,\dots,n$ are real-valued process, adapted to the filtration
of $x$. In the highly non-generic situation where the vector fields $\{X_i\}$ commute, $x_t$ becomes a function of 
$w_t$ and the problem trivializes. 
The argument in [1] sets up a duality between the processes $h$ and $r$, the lift of $\e$, in which
(in the non-commuting case)
the commutators $[X_i,X_j]$ play an explicit role.

The point of departure for the present work is the a priori selection of
an additional
 collection of
vector fields
$\{V_I: I\in {\cal I}\}$ on $M$ such that
 $$\{V_I(x) : I\in {\cal I}\} \k span \k T_xM, \k \forall x\in M.  \eqno (1.3)
$$
Thus in the elliptic case $\{V_I\}$ can be taken to be the
 set $\{X_1,\dots X_n\}$, whereas in the {\it hypoelliptic} case
(where the diffusion process (1.1) is degenerate but H\"ormander's condition
holds), one
can choose $\{V_I\}=
Lie(X_1,\dots,X_n)$, the Lie algebra generated by the vector fields $X_1,\dots,X_n$.
We construct  admissible  vector fields on $C_o(M)$ in the form
$$\eta_t = \sum_{I\in {\cal I}}h_I(t)V_I(x_t). $$
Somewhat
surprisingly,
it proves
to be possible to trade ellipticity in $\{X_1,\dots X_n\}$ for ellipticity
in $\{V_I\}$.
This enables us to establish our results under very general hypotheses.

The layout of the paper is as follows. Section 2 contains background material. The
results here are well-known, for the most part. Theorem 2.1 asserts that Riemann
integrals of continuous  adapted
paths have divergence given by an It\^o integral, while Theorem 2.2
states that It\^o  integrals with continuous
 adapted skew-symmetric  integrands are divergence free.
 The former result follows easily from the Girsanov theorem, the latter from the infinitesimal
rotation-invariance of the Wiener measure. Theorem 2.5 gives a relationship between
a  vector field $\eta$ along the path $x$ and  the lift of $\eta$ to the Wiener space. This
relationship, expressed in terms of the derivative of the stochastic flow of
 the SDE (1.1) and the inverse flow, plays a key role in Section 3.

Section 3 contains the main results of the paper. Theorem 3.1 gives the construction of a
class of vector fields on $C_o(M)$ as functions of $x$, under hypotheses that allow
the SDE (1.1) to be degenerate. The proof of Theorem 3.1 follows the above outline, and is an
extension of the
 argument in [1]. An essential step in the proof is the
decomposition of   non-tensorial terms in the lift obtained from Theorem 2.5,
into {\it tensorial} plus {\it skew-symmetric} parts.

 Theorem 3.2 is a variation on
Theorem
3.1 that exhibits a vector field on $C_o(M)$ with  given divergence. In particular,
we obtain a class of vector fields with divergence expressed in terms of
Ricci curvature. The interest of
 this result lies in the fact that formulae
of this type arise in the work of other authors,
 e.g Driver [6]  and Elworthy, Le Jan \& Li [8],
where they are obtained using  different methods.
In Example 3.3, Theorem 3.2 is applied to obtain vector fields with
divergence having no
extraneous dependence on the Wiener path $w$. This property is important
 in applications of the theorem
 that require
a degree of regularity  of the divergence
such as the study of quasi-invariance (this
point is discussed in the remark directly preceding Example 3.3).
 Theorem 3.4 is an intrinsic formulation of  Theorem 3.1 that does not depend
on the choice of a basis $\{V_I\}$. The proof of this result requires the 
introduction of a tensor that enables us to express the
Levi-Civita connection on $M$ in terms of a connection
 intrinsic to the diffusion process (1.1). In Theorem 3.7, we apply our theory to gradient 
systems. As a consequence (Corollary 3.8), we obtain Driver's result cited above.

In Section 4, we consider the special case where  the vector fields
$X_1,\dots,X_n$ are linearly independent. In this case, the problem under
consideration simplifies considerably and our argument simplifies accordingly. We 
conclude with an example where the SDE (1.1) takes values in the
Heisenberg group $G$. In this case we obtain explicit formulae for a class of
admissible vector fields $\eta$ on $C_o(G)$.

\section{Background material}
\n {\it  2.A Divergence theorems for Wiener space}

\n We present two such results. These concern the transformation of the Wiener measure
under Euclidean motions (the first under translations, the second under rotations).

 Let $\Omega$ denote the measure space for the Wiener process, equipped with the filtration
$${\cal F}_t = \sigma\{w_s/\k s \le t\}.
   $$

\begin{thm} Let $h:\Omega\times [0,T]\mapsto {\bf R}^n$ be a  continuous
adapted path.
Then the process $\int_0^\cdot h$ is admissible (with respect to the Wiener measure) and
$$
Div \Big [ \int_0^.h_sds\Big ]=\int_0^Th_s\cdot dw_s
$$
where $\cdot$  on the right of the equation denotes the Euclidean inner product.
\end{thm}

\n {\it Proof}. The result follows easily from the Girsanov theorem, which implies
that for $\Phi \in C_b^\infty\big (C_0({\bf R}^n)\big )$ and $\epsilon \in {\bf R}$,
$$
E\Big [\Phi(w+\epsilon\int_0^\cdot h_sds)\Big ]
=E[\Phi(w)G_\epsilon (w)]  \eqno (2.1)
$$
where
$$
G_\epsilon (w)\equiv \epsilon\int_0^Th_s\cdot dw_s -
 {\epsilon^2\over 2}\int_0^T||h_s||^2ds.
$$
Differentiating each side of (2.1) wrt $\epsilon$ and setting $\epsilon=0$ gives the
theorem.

\begin{thm} Let $A: \Omega\times [0,T]\mapsto so(n)$ be a  continuous adapted
process.
Then  $\int_0^\cdot Adw$ is admissible and
$$Div\Big [\int_0^.Adw\Big ] =0.
$$
\end{thm}

\n {\it Proof}. Define a process $\theta_t^\epsilon =\exp\epsilon (A_t)$ where $\exp$
denotes matrix exponentiation. Then $\theta_t^\epsilon$ is an
adapted $O(n)$-valued matrix  process with
$\theta_t^0=I. $  It follows from the infinitesimal rotation-invariance
of the Wiener measure that the law of the process
$$
w^\epsilon \equiv \int_0^\cdot \theta_t^\epsilon dw_t
$$
is invariant under $\epsilon$. Hence for $\Phi \in C_b^\infty(C_0){\bf R}^n))$,
we have
$$
E[\Phi(w^\epsilon)] = E[\Phi(w)].
$$
As before, differentiating in $\epsilon$ and setting $\epsilon=0$ gives the result.

\k

\k

{\n \it 2.B Geometric preliminaries}

\n In this section we introduce some geometric machinery that
will be needed in Section 3.  We adopt the summation convention throughout the paper:
whenever an index in a product (or
a bilinear form) is repeated, it will be assumed to be summed on.

First, let  $[g_{jk}]$ be the Riemannian metric defined on $M$ by
$$
g^{jk} =a_I^ja_I^k
$$
where
$$
V_I=a_I^j{\partial\over\partial x_j},\k I\in {\cal I}
$$
is the expression of the vector fields in local coordinates
(note that the matrix $[g^{jk}]$
 is non-degenerate by the spanning condition (1.3)).

 Denote the
corresponding
  inner product
 by $(.,.)$.
It is easy to see that
$$
V=(V,V_I)V_I, \k \forall V\in TM. \eqno (2.2)
$$
Let  $\tilde\nabla$ denote  the Levi-Civita
covariant derivative
corresponding
 to this metric.

The following constructions were introduced by Elworthy, Le Jan and Li (cf. [8]).
 Assume the set
 of vectors $\{X_1(x),\dots, X_n(x)\}$ span a subspace
$E_x$ of $T_xM$ of
 constant dimension as $x$
varies in $M$ and define $E$ to be the subbundle of $TM$
$$
E=\bigcup_{x\in M}E_x.
$$

Then $E$ becomes a Riemannian bundle under the inner product induced on $E$
 by the linear maps
$$
X(x): (h_1,\dots ,h_n) \in {\bf R}^n \mapsto
h_iX_i(x) \eqno (2.3)
$$
from the Euclidean space ${\bf R}^n$.

There is a  metric connection $\N$
on $E$ compatible with the metric $<.,.>$. This connection (termed the {\it Le 
Jan-Watanabe connection} in [8]),
is defined by
$$
\nabla_VZ =X(x)d_V(X^*Z), \kern 5 pt Z\in \Gamma(E), V\in T_xM,
$$
where $d$ represents the
derivative of the function
$$x\in M\mapsto X(x)^*Z(x)
\in {\bf R}^n.
$$

 The corresponding {\it Riemann curvature} tensor is defined by
$$
R(X,Y)Z = \N_X\N_YZ-\N_Y\N_XZ-\N_{[X,Y]}Z,
$$
and  the {\it Ricci tensor} by
$$
Ric(X) = R(X,e_i)e_i
$$
where $\{e_i\}$ is an orthonormal basis of $E_x$.

\begin{lem} \hfil\break
(i)\k  $<Y,X_i>X_i = Y,\k \forall Y \in E. \hfil\break
(ii) \k Ric(Y) = R(Y,X_i)X_i, \k \forall Y \in TM$.
\end{lem}

\n {\it Proof}. See [3. Sec. 2]

\k

\k

{\n \it 2.C Flow-related theorems}
\begin{lem} Let $g_t: M\mapsto M$ denote the stochastic flow  $x_0\mapsto x_t$
defined by the SDE (1.1).
Define
 $Y_t: T_{x_0}M \mapsto T_{x_t}M$ and $Z_t: T_{x_t}M \mapsto 
T_{x_0}M$ by
$
Y_t \equiv dg_t
$and $
Z_t\equiv Y_t^{-1}
$. Let $B$ denote a vector field on $M$ and $d$ the stochastic time
differential.
Then
$$
d\big [Z_tB(x_t)\big ] = Z_t\Big ([X_i,B](x_t)\circ dw_i+[V,B](x_t)dt\Big ).
$$
\end{lem}

\n {\it Proof}.

\k

\n Let $D_t$ denote the stochastic covariant differential along the path $x_t$,
with respect to the Levi-Civita $\tilde\N$ connection defined above. Then
 differentiating
with respect to the initial point $o$ in (1.10) gives\footnote{From this point on, we
 assume that all vector fields appearing in the equations are evaluated at
$x_t$. }
$$
D_tY= \tilde \N_{Y_t}X_i\circ dw_i +\tilde\N_{Y_t}Vdt.
$$
We then have
$$
D_tZ = D_t(Y_t^{-1})
$$
$$
=-Z_tD_tYZ_t
$$
$$
= -Z_t\Big (\tilde\N_{Id_t}X_i\circ dw_i + \tilde\N_{Id_t}Vdt\Big )
$$
where $Id_{t}$ denotes the identity map on $T_{x_t}M$. Thus
$$
d\big (Z_tB\big ) = D_tZB+Z_t\tilde\N_{dx_t}B
$$
$$
=-Z_t\Big (\tilde\N_{B}X_i\circ dw_i + \tilde\N_BVdt\Big ) +Z_t\Big
 (\tilde\N_{X_i}B\circ dw_i+\tilde\N_VBdt\Big )
$$
$$
d\big [Z_tB(x_t)\big ] = Z_t\Big ([X_i,B](x_t)\circ dw_i+[V,B](x_t)dt\Big ).
$$
as required.

\begin{thm} Let $r: \Omega\times [0,T]\mapsto {\bf R}^n$ be an It\^o
process.
Then the path $\eta\equiv dg(w)r$ is given by
$$
\eta_t = Y_t\int_0^tZ_sX_i(x_s)\circ dr_i \eqno(2.4)
$$
\end{thm}

\n {\it Proof}.  Note that $\eta$ is a vector field along the path $x$. Let
$U_s:T_o\mapsto
T_{x_s}M$ denote stochastic parallel translation along $x$.

Differentiating in (1.1)  with respect to $w$ gives
the following covariant
 equation for $\eta$
$$
D_t\eta=\tilde\N_\eta X_i(x_t)\circ dw_i+X_i(x_t)\circ dr_i
+\tilde\N_\eta V(x_t)dt \eqno (2.5)
$$
$$
\eta_0=0. \kern 180 pt
$$
We write (2.5) as
$$
d(U_t^{-1}\eta) = U_t^{-1} \tilde\N_\eta X_i(x_t)\circ
dw_i+U_t^{-1}X_i(x_t)\circ dr_i
+U_t^{-1}\tilde\N_\eta V(x_t)dt.
$$
Denoting the path $t\mapsto U_t^{-1}\eta_t $ by $y$, we note that the equation for
$y$ has the form
$$
dy= M_i(t)y_t\circ dw_i +M_0(t)y_t + U_t^{-1}X_i(x_t)\circ dr_i \eqno (2.6)
$$
where $M_j(t), j=1,\dots, n$ are linear operators on $T_{o}M$.

On the other hand, differentiation in (1.1) with respect the the initial point $o$
gives
 the following equation for $\tilde Y_t\equiv U_t^{-1}Y_t  $
$$
d\tilde Y =M_i(t)\tilde Y_t\circ dw_i+M_0(t)\tilde Y_tdt \eqno  (2.7)
$$
$$
\tilde Y_0=I.\kern 110 pt
$$
Equation (2.6) can be solved in terms of $\tilde Y$ using an operator version
of the familiar
``integrating factor" method used to solve first order linear ODE's.
Noting, then, that $\tilde Y^{-1}$ is an integrating factor for (2.6) and using this
to solve for $y$
gives
$$
y_t =\tilde Y_t\int_0^t \tilde Y_s^{-1}U_s^{-1}X_i(x_s)\circ dr_i. \eqno (2.8)
$$
Writing (2.8) in terms of $\eta$ and $Y$, we obtain (2.4).

\k

\n {\bf Remarks}

\n 1. Theorem 2.5 gives an alternative proof of the ``lifting" equation (3.2)
in [1].

\k

\n 2. Suppose $\eta$ in (2.4) has the form $\eta_t=X_i(x_t)h_i(t)$ for an
 ${\bf R}^n$-valued process $h = (h_1,\dots, h_n)$. Then, writing
$$X=[X_1\dots X_n]$$
and solving for $dr$ in (2.4), we have
$$
Z_tX(x_t)\circ dr=d\big [ Z_tX(x_t)h_t\big ].
$$
This equation suggests that  $r$ can be considered as a type of ``covariant derivative"
of $h$ along $x$, where the operator $Z_tX(x_t)$ plays the role of backward parallel
translation.

\section{Divergence theorems for degenerate diffusions}

\n {\it 3.A A divergence theorem}

\n Let $X$ be as  defined in (2.3). Then the SDE (1.1) may be written
$$
dx = X(x_t)\circ d\tilde w
$$
where
$$
d\tilde w=dw+X(x_t)^*V(x_t)dt
$$
and  the adjoint map is defined using the metric $<.,.>$ on $E$ (so $X(x)^*$
is a left inverse
for $X(x)$).
By the Girsanov theorem , the law $\tilde \nu$ of  
of $\tilde w$ is equivalent to the law $\nu$ of $w$, with Radon-Nikodym derivative
${d\tilde \nu \over d\nu}$ given by
$$
G(w) =\exp\big (\int_0^T X(x_t)^*V(x_t)\cdot dw -
 {1\over 2}\int_0^T||X(x_t)^*V(x_t)||^2dt\Big ).
$$
Suppose that  $r$ is an admissible lift for the vector field $\eta$
under the map
 $\tilde g: \tilde w \mapsto x$. Then 
$$
E[\eta\phi(x)] = E\big [G(w)\cdot r(\Phi\circ\tilde g)(w)\big ]
$$
$$
=E\big [\Phi\circ\tilde g(w)Div (G\cdot r)\big ].
$$
$$
E\big [\Phi\circ \tilde g(w)\{G\cdot Div(r) - r(G)\}].
$$
Thus  $\eta$ is admissible.

In view of this discussion, there is no loss in generality in assuming $V=0$ and we
shall assume in the sequel that this is the case.\footnote{It is clear that our argument
 will work for non-zero drift, however this reduction to the case $V=0$
simplifies the later calculations.}

\k

We introduce the following  tensors 
$\{S_I\}$ and $\{T_I\}$  associated to the
vector fields
$\{V_I\}$
$$
S_I(X)= \N_{V_I}X +[X,V_I],\k X\in E,
$$
and
$$
T_I(X)=S_I(X) -<\N_{V_I}X_i,X>X_i,\k X \in E. 
$$
\begin{thm} Let $r=(r_1, \dots, r_n)$ be a path
in the Cameron-Martin space of ${\bf
R}^n$ and define
$\{h_I: I\in {\cal I}\}$ by the linear stochastic  system
$$
dh_I=(X_i,V_I)  \dot r_idt - (T_J (\circ dx),V_I) h_J \eqno (3.1)$$
$$
\kern 40 pt h_I(0) =0. \kern 148 pt
$$
Then the vector field $\eta_t\equiv h_I(x_t)h_I(t), \kern 3 pt t\in [0,T]$
 is admissible on $C_o(M)$.
\end{thm}

\n {\it Proof}.

\k

\n We first note that Theorem 2.5 implies that $r$ is lift of $\eta$ if $r$
satisfies
$$
X_i dr_i =Y_t d\big [Z_t\eta_t\big ] \eqno(3.2)
$$
Substituting $\eta_t=h_I(t)V_I(x_t)$ into (3.2) and using Lemma 2.4, we
have
$$
X_i dr_i=V_I\circ dh_I + [X_j, V_I]h_I\circ dw_j \eqno (3.3)
$$
Writing the Lie bracket term involving $X_j$ in terms of the connection
$\nabla$ and using Lemma 2.3 (i) gives
$$
[X_j,V_I]=S_I(X_j)-\N_{V_I}X_j
$$
$$
=S_I(X_j) -<\N_{V_I} X_j,X_i>X_i
$$
Denote
$$
G_I^{ij} = <\N _{V_I}X_i,X_j> -<\N _{V_I}X_j,X_i> \eqno (3.4)
$$
Then we have
$$
[X_j,V_I]  = G_I^{ij}X_i + T_I(X_j)
$$
Substituting this into (3.3) gives
$$
X_i dr_i=V_I\circ dh_I +G_I^{ij}h_IX_i\circ dw_j +T_I(\circ dx)h_I. \eqno (3.5)
$$
We note that, more generally, a semimartingale path $\tilde r$ is a lift of $h_IV_I$ if 
equation (3.5) holds with the left hand side replaced by the Stratonovich 
differential $X_i\circ d\tilde r_i$.

Suppose now the coefficient functions $\{h_I\}$ satisfy the system
$$
X_i dr_i=V_I\circ dh_I+T_I(\circ dx)h_I, \eqno (3.6)
$$

\kern -15 pt
$$
h_I(0)=0.  \kern 75 pt
$$
Then
$$
X_i\big [  dr_i +G_I^{ij}\circ h_I\circ dwj\big ]=V_I\circ dh_I+
G_I^{ij}X_ih_I\circ dw_j+
T_I(X_j)h_I\circ dw_j
$$
So if we  define
$$
\tilde r_i = r_i+ \int_0^\cdot G_I^{ij}h_I\circ dw_j. \eqno (3.7)
$$
then (3.3) holds with $ r$ replaced by $\tilde r$. It follows that  $\tilde r$ is a
 lift of $\eta$,
where
$$\eta_t = h_I(t)V_I(x_t). \eqno (3.8)$$
Furthermore, the the skew-symmetry of the
 functions $G_I^{ij}$ in the
upper indices and Theorem 2.2 imply that the Stratonovich integral in (3.7)
can be written as  a  Riemann integral  plus a
divergence-free It\^o integral. It follows from Theorems 2.1 and 2.2 that $\tilde r$ is admissible.
Note also  that by (2.2), the processes $h_I$ defined by (3.1)
 satisfy equation (3.3).

We have thus shown that $\tilde r$ is an admissible lift
 to the Wiener space of the vector field
$
\eta
$
 in (3.8). In view of  Definition 1.2, we have  for any test function $\Phi$ on $C_o(M)$
$$
E\big [(\e\Phi)(x)\big ] = E\big [r(\Phi\circ g)(w)\big]$$
$$ = E\big [\Phi \circ
g(w)
 Div (r)\big ]
$$
$$
= E\Big [ \Phi(x)E\big [Div(r)/x\big ]\Big ].
$$
Thus $\eta$ is admissible and
$$
Div(\e)(x) = E\big [Div(r)/x\big ].
$$

\k

\n  {\it 3.B Computation of the divergence}

\n In order to compute the divergence of the vector field $\eta$ in Theorem 3.1,
it is necessary to
 convert the Stratonovich
integral in (3.7) into It\^o form. The relation between the
Stratonovich and It\^o differentials is formally
$$
G_I^{ij}h_I\circ dw_j = G_I^{ij}h_I dw_j + {1\over
2}d(G_I^{ij}h_I)dw_j. \eqno (3.9)
$$

Write
$$\alpha_I^{kij} = <\N _{X_k}\N _{V_I}X_i,X_j>+ <\N_{V_I}X_i,\N_{X_k}X_j>
$$
$$
-<\N _{X_k}\N _{V_I}X_j,X_i>
 - <\N_{V_I}X_j,\N_{X_k}X_i> \eqno (3.10)
$$
and
$$
\beta_I^k = -(T_J(X_k),V_I)h_J. \eqno (3.11)
$$
Then by (3.1) and (3.4)
$$
dG_I^{ij}=\alpha_I^
{kij}dw_k+ \{\dots\}d
$$
and
$$
dh_I=\beta_I^kdw_k + \{\dots\}dt.
$$
Substituting these into (3.9) and using the It\^o rules
$$dw_idw_j=\delta_{ij}dt,\k dw_idt=0$$ we see
that the Ito-Stratonovich correction term in (3.9) is
$$
{1\over 2}(\alpha_I^{kik}h_I+G_I^{ik}\beta_I^k)dt. \eqno (3.12)
$$
Thus  (3.7) becomes
$$
\tilde r_i = r_i +\int_0^\cdot G_I^{ij}h_Idw_j +{1\over 2}\int_0^\cdot
(\alpha_I^{kik}h_I+G_I^{ik}\beta_I^k ) dt.
$$
As remarked in the proof of Theorem 3.1, the   It\^o integral has divergence
 zero and using
Theorem 2.1 we obtain
$$
Div(\tilde r) = \int_0^T \Big (\dot r_i+{1\over 2}\big
(\alpha_I^{kik}h_I+G_I^{ik}\beta_I^k\big )\Big )dw_i
$$
Hence
$$
Div(\eta)= E\Big[\int_0^T \Big (\dot r_i+
{1\over 2}\big  (\alpha_I^{kik}h_I+G_I^{ik}\beta_I^k\big )\Big )dw_i\big /x\Big ] \eqno (3.13)
$$
where the $\alpha$'s and $\beta$'s are given in (3.10) and (3.11).

\k

\k

By adjusting the right hand side in equation (3.1) by the addition of a suitably
chosen drift term, the above argument
can
easily be
modified to give

\begin{thm}
Let $\gamma: \Omega\times [0,T]\mapsto {\bf R}^n$ be a $C^1$ adapted process
and define $\{h_I\}$ by $h_I(0)=0$ and
$$
dh_I= \Big (\big (d\gamma_i-{1\over 2}G_J^{ik}\beta_J^kdt\big )X_i
+\big (T_J(\circ dx)-{1\over 2}\alpha_J^{kik}X_idt\big )h_J,V_I\Big ).
$$
Then the vector field $\eta_I = h_IV_I$ is admissible and for every test function $\Phi$
on $C_o(M)$, we have
$$
E\big [  (Z\Phi)(x)\big ]=E\Big [\Phi(x)\int_0^T \dot \gamma_idw_i\Big ].
$$
\end{thm}

\n The proof of Theorem 3.2 is an easy modification of the argument above, where we
 replace
 $r$ by the path
$$
\tilde r_i  =  \gamma_i - {1\over 2}\int_0^\cdot
\big  (\alpha_I^{kik}h_I+G_I^{ik}\beta_I^k\big )dt.
$$
The essential point is that the correction term (3.12) in the computation of
 the divergence does not{\it  explicitly}
involve the path $r$.

\k

\n {\bf Corollary to Theorem 3.2} \k
{\it Given any path $r$ in the Cameron-Martin space of ${\bf R}^n$, we can construct an
admissible vector field $\eta$ on $C_o(M)$
such that  
$$
E\big [(\eta\Phi)(x)\big ] = E\Big [\Phi (x)\int_0^T\Big (\dot r_i +
 {1\over 2}<Ric(\eta),X_i>(x_t)\Big )dw_i\Big ]. \eqno (3.15)
$$
}

\n {\bf Remarks}

\n 1. Formula (3.12) is similar to those
 appearing in the work of Driver
 [6], [7] and Elworthy, Le Jan \& Li [8].

\k

\n 2. The appearance of the conditional expectation in (3.14) and (3.15)
entails a loss of
 information concerning the regularity of the function $Div(\e)$.
This point is crucial in certain
 applications of the results presented here.
For example, the  regularity of $Div(\eta)$   plays a major role
in recent work of the author [4]
in which the admissibility of $\eta$ is used, in the elliptic setting, to establish quasi-invariance of the
law of $x$ under the flow generated by $\eta$ on $C_o(M)$.

With this in mind, we note that by choosing the process $\gamma$  in (3.14)
appropriately, we can eliminate the {\it extraneous} dependence of the integral on $w$
and thus circumvent this problem. The next
example
 illustrates this point.

\k

\k

\n {\bf  Example 3.3} \k
 Suppose $B$ is a smooth vector field on $M, \rho$ is a   
deterministic $C^1$ real-valued function,  and define
$$\gamma_i(t)=\int_0^\cdot \rho_t(B, X_i)(x_t)dt $$
so
$$
\int_0^T\dot\gamma_idw_i=\int_0^T \rho_t(B, X_i)dw_i.
$$
Using the Levi-Civita connection $\tilde\N$ to write this in 
Stratonovich form we have
$$
\int_0^T \rho_t(B, X_i)dw_i =
$$
$$
\int_0^T \rho_t(B, X_i)\circ dw_i -{1\over 2}
\int_0^T\rho_t\Big ( \big (\tilde\N_{X_i}B, X_i\big )+\big (B,\tilde 
\N_{X_i}X_i\big )\Big )dt =
$$

$$
\int_0^T \rho_t(B, \circ dx) -{1\over 2}
\int_0^T\rho_t\Big ( \big (\tilde\N_{X_i}B, X_i\big )+\big (B,\tilde 
\N_{X_i}X_i\big )\Big )dt \eqno (3.16)
$$
Since  (3.16) is measurable
 with respect to $x$, (3.14) becomes
$$
Div(\eta) =
\int_0^T \rho_t(B, \circ dx) -{1\over 2}
\int_0^T\rho_t\Big ( \big (\tilde\N_{X_i}B, X_i\big )+\big (B,\tilde 
\N_{X_i}X_i\big )\Big )dt. 
$$
In particular, $Div(\eta)$ is an explicit function of the path $x$.

\k

\k

\n {\it 3.C A basis-free formulation of the argument}

\n Assume now that $M$ is a {\it Riemannian} manifold.  In this case we can
formulate the preceding argument {\it intrinsically}, i.e.
 in a way that does not depend
on the choice of a basis $\{V_I\}$. 

Let $\tilde \nabla$ denote the Levi-Civita covariant derivative with respect to the
Riemannian metric on $M$ and $\tilde D$ the corresponding covariant stochastic 
differential. As before, $<.,.>$ and  $\nabla$ will denote the inner product
and the connection on the subbundle $E$ 
introduced in Section 2.B. 

We define
$$
T(X,Y)=\tilde\N_YX- \N_YX,\k Y\in TM,X\in E, \eqno (3.17)
$$
noting that $T$ is tensorial in {\it both} arguments. 

Let $r: [0,T]\times \Omega\mapsto {\bf R}^n $ be an It\^o semimartingale 
$$
dr_k(t)=b^{kj}(t)dw_j+c^k(t)dt  
$$
where $b^{kj}$ and $c^k$ are adapted continuous processes. Then
differentiation in equation (1.1) gives
the following covariant equation for the path $\eta\equiv dg(w)r$ 
$$
\tilde D_t\eta = \tilde\N _\eta X_i\circ dw_i+X_i\circ dr_i
$$
$$
=\N_\e X_i\circ dw_i+T(X_i,\e )\circ dw_i +X_i\circ dr_i
$$
$$
=<\N_\e X_i,  X_j>X_j\circ dw_i+T(X_i,\e)\circ dw_i+X_i\circ dr_i
$$
$$
=  <\N_\e X_j,X_i>X_j\circ dw_i +G_\e^{ij}X_j\circ dw_i+T(X_i,\e )\circ dw_i+X_i\circ dr_i
$$
where
$$
G_V^{ij}\equiv <\N_V X_i,X_j>-<\N_V X_j,X_i>.
$$
Thus 
$$
\tilde D_t\e= <\N_\e X_j,X_i>X_j\circ dw_i +T(X_i,\e)\circ dw_i
+ X_i\big (\circ dr_i+ G_\e^{ji}\circ dw_j\big ). \eqno (3.18)
$$
We now have

\k

\n {\bf Theorem 3.4}\k {\it Let $r$ be any Cameron-Martin path in ${\bf R}^n$ and define a vector
 field $\e$ along
$x$ by the covariant SDE
$$
\tilde D_t\e= \big [<\N_\e X_j, \cdot >X_j +T(\cdot ,\e)\big ] (\circ dx ) 
+ X_i \dot r_i dt \eqno (3.19)
$$
$$
\e (0)=0. \kern 160 pt
$$
Then $\e$ is admissible and for any test function $\Phi$ on $C_o(M)$,
$$
E [(\e \Phi)(x)]=E\Big [\Phi(x)\int_0^T(\dot r_i-{1\over 2}\alpha_i)dw_i\Big ], \eqno (3.20)
$$
where
$$
\alpha_i(t) = \big [<\N_{X_k}(\N_\e X_k),X_i> + <\N_\e X_k, \N_{X_k}X_i>
$$
$$
-<\N_{X_k}(\N_\e X_i),X_k>-<\N_\e X_i, \N_{X_k}X_k>\big ] (x_t).
$$
}

\n {\it Proof}.  We note that equation (3.18)
implies  $\tilde r$ is a lift of $\eta$, where
$$
\tilde r_i =r_i-\int_0^\cdot G^{ji}_\eta \circ dw_j. \eqno (3.21)  
$$
Since the functions $G^{ji}_\e$  are skew-symmetric  in the indices
 $j$ and $i$,
 Theorems 2.1 and 2.2 imply that $\tilde r$ is an admissible vector
 field on the Wiener space. 
As before, for any test function $\Phi$
on $C_o()M)$,  we 
have
$$
E\big [ D\Phi (x)\eta\big ]= E\big [\Phi(x)Div(\tilde r)\big ].
$$
and it follows that $\e$ is admissible as claimed.

As in Section 2.B, the divergence $Div(\tilde r)$  is computed by converting 
the Stratonovich integrals in (3.21) into
It\^o form and applying Theorem 2.1. This yields (3.20) and so completes the proof.

\k

\k\n {\bf Remark 3.5 }

 \n It is clear that the argument used to
prove Theorem 3.4 is valid in more generality,  with the deterministic
 Cameron-Martin
 path $r$ replaced by any
 ($x$-measurable) random path of the form 
$$
 r =\int_0^\cdot A(s)dw_s+\int_0^\cdot B(s)ds. \eqno (3.22)
$$
where $A: \Omega \times [0,T] \mapsto so(n)$ and $B: \Omega \times [0,T] 
\mapsto {\bf R}^n
$  are continuous adapted processes. In view of Theorems 2.1 and 2.2, it is natural to consider
 the Wiener space
 $C_0({\bf R}^n)$ as a manifold with tangent  bundle 
$ \cup_w T_w C_0({\bf R}^n)$, where each fiber $T_w C_0({\bf R}^n)$
consists of paths of the form (3.22). 

For each such path $r =r(x)$,
equation (3.19) produces a vector field  $\eta$ on  $C_o(M)$ that is then lifted 
to a vector field  
 $\tilde r$ on $C_0({\bf R}^n)$ by equation (3.21). 
We summarize these constructions as follows. 

Define  
$$
H(r)= (r,\e), \k r\in TC_0({\bf R}^n)
$$
and let $$\pi: TC_0({\bf R}^n) \mapsto C_0({\bf R}^n)$$ denote the
bundle projection.

 Then the  chain of maps in  Theorem 3.4 and its proof  is illustrated by the 
 commutative diagram

\kern 5 pt

$$\xymatrix{ &TC_0({\bf R}^n)\times TC_o(M)\ar[rd]^{(3.21)} \\
TC_0({\bf R}^n)\ar[r]^{(3.19)}\ar[ur]^H \ar[d]_{\pi} &TC_o(M)  
 &TC_0({\bf R}^n)\ar[l]_{dg}\\
C_0({\bf R}^n)\ar[r]_g   &C_o(M)\ar[ul]^r\ar[u]^\eta  &C_0({\bf R}^n)\ar[l]^g\ar[u]_{\tilde r}
}
$$

\n {\it 3.D Gradient systems }

\n  Suppose $M$ is an isometrically embedded submanifold\footnote{By 
Nash's
embedding theorem, every finite-dimensional Riemannian manifold can be realized this way.} of a
Euclidean space ${\bf R}^N$. Define $X_i= Pe_i, \k 1\le i\le N$ where $e_1,\dots, e_N$ is
the standard orthonormal basis of ${\bf R}^N$ and $P(x)$ is  orthogonal
projection onto the tangent space $T_xM$. Then the diffusion process
$x$ in equation (1.1) is a Brownian motion in $M$.

In this case the connection $\N$ coincides  with the Levi-Civita 
connection on $M$ (cf.
[8]), hence the tensor $T$ defined in (3.17) is zero.
Equation (3.19) thus becomes
$$
\tilde D_t\e= <\N_\e X_j, \circ dx >X_j 
+ X_i \dot r_i dt \eqno (3.23)
$$
A further reduction results from

\k

\n{\bf  Lemma 3.6}\k {\it For all $V\in TM$ and $W\in {\bf R}^N$}
$$
<\N _VX_j,W>X_j = 0.
$$
\n {\it Proof}. Using the classical representation of
 the Levi-Civita connection and denoting the Frechet derivative by $d$, we have
$$<\N _VX_j,W>X_j = Pe_j<PdP(V)e_j ,W>$$
$$
=Pe_j<e_j,dP(V)PW>
$$
$$
=PdP(V)PW. 
$$
Differentiating the relation $P^2=P$ gives
$$
dP(V)P+PdP(V)=dP(V).
$$
Thus
$$
dP(V)P = dP(V)-PdP(V) = QdP(V)
$$
where $Q=I-P$. Then
$$
PdP(V)P = PQdP(V)=0
$$
and the result follows.

\k

In view of Lemma 3.6, equation (3.22) reduces to
$$
\tilde D_t\e = X_i dr_i.
$$
Hence
$$
\e_t=U_t\int_0^tU_s^{-1}X_i dr_i. \eqno (3.24)
$$
where $U$ denotes parallel translation along $x$.
This yields

\k

\n {\bf Theorem 3.7} \k {\it If $r$ is any (random, $x$-adapted) path such that $\dot r \in L^2[0,T]$ then the vector field $\e$ defined by 
(3.24) is admissible. }

\k

 In particular, let $h$ be any path in the Cameron-Martin space of $T_o(M)$ 
 and define
$$
r_i=\int_0^\cdot <U_t\dot h_t,X_i>dt,\k i=1,\dots, N.
$$
Then the integral in (3.24) becomes $h_t$ and we obtain the following  result of
 Driver (cf. [6])

\k

\n {\bf Corollary 3.8}\k {\it For every path $h$ in the Cameron-Martin space of $T_o(M)$, 
the vector field
$\e_t \equiv U_th_t$ is admissible.}

\k

Finally, we note that every adapted vector field on $C_o(M)$ with an admissible
lift to the Wiener space  is obtained from Theorem 3.4.  Denote
the process $\e$ in Theorem 3.4 by $\e^r$. Then we have

\k

\n{\bf  Proposition 3.9}\k  {\it Suppose $\e$ is an adapted vector field
 on $C_o(M)$ 
such that
$$
\e = dg(w)r
$$
for some $r\in TC_0({\bf R}^n)$. Then there exists $\bar r \in  TC_0({\bf R}^n)$ such that
$\e = \e^{\bar r}$.}

\k

\n {\it Proof}. This follows immediately from equations (3.18) and (3.19).
 We define
$\tilde r$ by
$$
\tilde r_i = r_i + \int _0^\cdot G_\e ^{ji} \circ dw_j,\k i = 1,\dots, n.
$$

\section{ Linearly independent diffusion coefficients}

\n In this section we consider the special case where the vectors
$\{X_1(x),\dots,X_n(x)\}$ are linearly
 independent at every point $x \in M$\footnote{In the {\it elliptic} 
case there is a topological obstruction to this condition,
i.e. if $M$ has non-zero Euler characteristic then it is impossible. However,
 the condition is reasonable in the non-elliptic case.}. 
As we shall see, this implies that the Wiener path $w$ is a function of
the solution $x$ of the SDE (1.1) i.e.
$$
w=\Theta(x)
$$
where $\Theta$ is a measurable function on $C_o(M)$. In this case
 the
following simplified version of the method used in Section 3 
produces admissible vector fields on $C_o(M)$. 

Choose
$r$ to be any process of the form 
$$
r_t=\int_0^tA(s)dw_s+\int_o^tB(t)dt, \k t\in [0,T] \eqno (4.1)
$$
where $A$ and $B$ are continuous adapted processes with values in $so(n)$ and 
${\bf R}^n$
and define $\eta$ by (2.4), i.e.
$$
\eta_t=Y_t\int_0^tZ_sX_i(x_s)\circ dr_i.
$$
By Theorems 2.1, 2.2 and 2.5, $r$ is an admissible lift of 
$\eta$, hence $\eta(w)=\eta(\Theta(x))$ is an admissible vector field on
$C_o(M)$. 

We now study how the formulae in Section 3 reduce in the linearly independent
case. As before, define $X(x): {\bf R}^n \mapsto T_xM$ by 
$$
X(x)(h_1,\dots, h_n) =X_i(x)h_i.
$$
 We will need the following result.
\begin{lem} The vectors $X_1(x),\dots, X_n(x)$ are linearly independent
if and only if
$$
X(x)^*X(x)=I_{{\bf R}^n}.
$$
\end{lem} 

\n Since Lemma 4.1 is elementary, the proof will be omitted.  

\k

Assume now  that
$\{X_1,\dots,X_n\}$ are linearly independent. 
Then Lemma 4.1 enables us to solve the SDE (1.1) for $w$ in terms of $x$
 and obtain
$$
dw=X(x_t)^*\circ dx,
$$
thus $w=\theta(x)$, as claimed above. We also have

\k

\n {\bf Corollary to Lemma 4.1}\k  {\it For $a_i\in C^\infty (M), i = 1,\dots, n$
and $V\in TM$
$$
\nabla_V(a_iX_i)=V(a_i)X_i.
$$
In particular
$$
\nabla_VX_i=0, \k i=1,\dots, n.
$$}

 The corollary implies that the functions $G_I^{ij}$ in (3.4) are all zero.
Furthermore, the tensors $T_I$ in Section 3 take the form
$$
T_I(aX_i) = a[X_i,V_I],\k i = 1,\dots, n
$$
for $a\in C^\infty (M)$.
Theorem 3.1 then becomes
\begin{thm}  Suppose the process $r$ is defined as in (4.1) and the
functions $h_I$ are chosen to satisfy
$$
dh_I=(X_i,V_I)\circ dr_i-\big ([X_i,V_J],V _I\big )h_J\circ dw_i \eqno (4.2)
$$
$$
  h_I(0)=0. \kern 140 pt
$$
Then the vector field $\eta=h_IV_I$ is admissible and
$$
Div(\eta) = \int_0^TB_i(t)dw_i.
$$
\end{thm}

\k

\n {\bf Example 4.3} 

\n Let $M$ be the {\it Heisenberg group},
 i.e. the Lie group
${\bf R}^3$ with group multiplication 
$$
 (a_1,a_2,a_3)\cdot (b_1,b_2,b_3)=\Big (a_1+b_1, a_2+b_2,a_3+b_3+
{1\over 2}(a_1b_2-b_1a_2)\Big ).
$$
Let
$$
X_1={\partial\over \partial x}-{y\over 2}{\partial \over\partial z}
$$
$$
X_2={\partial\over \partial y}-{x\over 2}{\partial \over\partial z}
$$
 and define $V_1 = X_1, V_2=X_2$, and 
$$V_3 = [V_1,V_2] = {
\partial\over 
\partial z}.
$$
Then 
$$\kern -5 pt[X_1,V_2]=V_3$$
$$ [X_2,V_1]=-V_3$$
$$\kern 35 pt[X_i, V_j]= 0,\k i+j\ne 3. $$
Thus equation (4,2), which we write in the form
$$
V_I\circ dh_I=X_i\circ dr_i-[X_i,V_I]h_I\circ dw_i
$$
becomes
$$ V_1\circ dh_1+V_2\circ dh_2+V_3\circ dh_3$$
$$
 =X_1\circ dr_1 +X_2 \circ dr_2 +V_3(h_1\circ dw_2-h_2\circ dw_1). \eqno (4.3)
$$
Since the vectors  $\{V_1,V_2,V_3\}$ are linearly independent
 equation (4.3) has 
a unique solution, given by
$$
  \kern - 93 pt          h_1=r_1
$$
$$
  \kern - 93 pt           h_2 =r_2
$$
$$
           h_3 = \int_0^\cdot r_1\circ dw_2 -r_2\circ dw_1. \eqno (4.4)
$$
As point of interest, we note that if  $(w_1,w_2)$ is substituted for $(r_1,r_2)$ 
then the integral in (4.4) becomes the
{\it Levy area} (this  is not, however, 
an admissible choice of $r$).

\vfil\break

\centerline {\bf References}

\k

\n [1] D. Bell,      Divergence theorems in path space. {\it J. 
Funct. Anal}.
 218 (2005), no. 1, 130-149.

\kern 3 pt

\n [2] D. Bell, Divergence theorems in path space II: degenerate diffusions,
{\it C. R. Acad. Sci. Paris Sr. I}  342 (2006),  869-872.

\kern 3 pt

\n [3] D. Bell, {\it The Malliavin Calculus}, 2nd edition. Dover Publications, Mineola,
NY, 2006.

\kern 3 pt

\n [4] D. Bell, Quasi-invariant measures on the path space of a diffusion,
{\it C. R. Acad. Sci. Paris Sr. I} 343 (2006), 197-200.

\kern 3 pt

\n  [5] J.-M.Bismut, {\it Large deviations
and the Malliavin calculus}. Progress
in Mathematics, 45. BirkhŠuser Boston, Inc., Boston, MA, 1984.

\kern 3 pt

\n  [6] B. Driver, A Cameron-Martin type quasi-invariance
theorem for Brownian motion on a compact manifold. {\it
J. Funct. Anal.}  109 (1992) 272-376.

\kern 3 pt

\n [7] B. Driver, A Cameron-Martin type quasi-invariance theorem
 for pinned Brownian motion on a compact Riemannian manifold (1994).

\kern 3 pt

\n [8] K. D Elworthy, Y. Le Jan and X.-M. Li, On
the Geometry of Diffusion Operators and Stochastic
Flows, Lecture Notes in Mathematics  1720,
Springer-Verlag, 1999.

\kern 3 pt

\n [9] O. Enchev and D. Stroock, Towards a Riemannian geometry on the path space over a
 Riemannian manifold. 
{\it J. Funct. Anal. } 134  (1995),  no. 2, 392--416.

\kern 3 pt

\n [10] E. P. Hsu,  Quasi-invariance of the Wiener measure on the path
space over a
compact
 Riemannian manifold.  {\it
J. Funct. Anal.} 134 (1995) 417-450.

\kern 3 pt

\n [11] E. P. Hsu, Quasi-invariance of the Wiener measure on path spaces: noncompact case.
 J. Funct. Anal.  193  (2002),  no. 2, 278--290.

\kern 3 pt

\n [12] P. Malliavin, Stochastic calculus of variations and
hypoelliptic operators. {\it Proceedings of the International Conference on
Stochastic Differential Equations}, Kyoto, 195-263. 
Kinokuniya
and Wiley, 1976.

\kern 3 pt

\n [13] P.Malliavin,
{\it Stochastic Analysis}.
Grundlehren der Mathematischen Wissenschaften
 [Fundamental Principles of Mathematical Sciences], 313. Springer-Verlag, Berlin, 1997.

\kern 3 pt

\n [14] D. Nualart, {\it The Malliavin calculus and related topics}.
Probability and its Applications (New York). Springer-Verlag, New York, 1995.

\kern 3 pt

\n [15] S. Watanabe, Lectures on stochastic differential equations and Malliavin calculus.
 Notes by M. Gopalan Nair and B. Rajeev. Tata Institute of Fundamental
Research Lectures on Mathematics and Physics, 73.   Springer-Verlag, Berlin, 1984.

\end{document}